\documentclass[11pt]{amsart}

\usepackage{latexsym}
\usepackage{graphicx}
\usepackage{psfrag}
\usepackage{amssymb}

\newcommand{\comment}[1]{}%$\star$\textsc{#1}$\star$}

  \newtheorem{thm}{Theorem}[section]
  \newtheorem{defn}[thm]{Definition}
  \newtheorem{lemma}[thm]{Lemma}
  \newtheorem{cor}[thm]{Corollary}
  \newtheorem{prop}[thm]{Proposition}

  \newtheorem{example}[thm]{Example}
  \newtheorem{remark}[thm]{Remark}

\def\proof{{\noindent{\it Proof.\ }}}
\def\endproof{{\hfill $\Box$}}
\def\implies{\Rightarrow}

\def\<{\langle}
\def\>{\rangle}
\def\0{{{\bf 0}}}
\def\CE{{\mathcal E}}
\def\CH{{\mathcal H}}
\def\EE{{\underline E}{\scriptstyle \,}}
\def\HH{{\widetilde H}{}{}}

\def\NN{{\mathbb N}}
\def\RR{{\mathbb R}}
\def\ZZ{{\mathbb Z}}

\def\SS{{\mathcal S}}
\def\UU{{\mathcal U}}
\def\gp{{\rm gp}}
\def\mm{{\mathfrak m}}
\def\oF{{\overline F}{}}
\def\oQ{{\overline Q}{}}
\def\ol#1{{\overline {#1}}}
\def\pp{{\mathfrak p}}

\def\wt#1{{\widetilde{#1}}}
\def\xx{{\mathbf x}}
\def\yy{{\mathbf y}}
\def\ext{{\rm Ext}}
\def\hom{{\rm Hom}}
\def\qgp{{Q^{\rm gp}}}
\def\sat{{\rm sat}}

\def\too{\longrightarrow}
\def\van{{\operatorname{van}}}
\def\ModR{{\mathcal M}}
\def\cech{{\check{C}}}
\def\eext{{\operatorname{\underline{Ext}}}}
\def\hhom{{\operatorname{\underline{Hom}}}}
\def\spot{{\hbox{\raisebox{1.7pt}{\large\bf .}}\hspace{-.5pt}}}
\def\coker{{\rm coker}}
\def\nothing{\varnothing}

	{\begin{list}
		{\noindent\makebox[0mm][r]{\arabic{enumi}.}}
		{\usecounter{enumi} \topsep=1.5mm \itemsep=0mm}
	}
	{\end{list}}

	{\begin{list}
		{\noindent\makebox[0mm][r]{\arabic{enumi}.}}
		{\usecounter{enumi} \topsep=1.5mm \itemsep=-1mm}
	}
	{\end{list}}

\begin{document}%%%%%%%%%%%%%%%%%%%%%%%%%%%%%%%%%%%%%%%%%%%%%%%%%%%%%%%%%
%%%%%%%%%%%%%%%%%%%%%%%%%%%%%%%%%%%%%%%%%%%%%%%%%%%%%%%%%%%%%%%%%%%%%%%%%

\title{Bass Numbers of Semigroup-Graded Local Cohomology\\
	{\footnotesize \rm \textsc{(15 October 2000)}}}
\author{David Helm}
\author{Ezra Miller}
\address{University of California at Berkeley, Massachusetts Institute of
Technology}
\email{dhelm@math.berkeley.edu, ezra@math.mit.edu}
%\keywords{}
%\subjclass{}

\begin{abstract}
\noindent
Given a module $M$ over a ring $R$ which has a grading by a semigroup
$Q$, we present a spectral sequence that computes the local cohomology
$H^i_I(M)$ at any graded ideal $I$ in terms of $\ext$ modules.  This
method is used to obtain finiteness results for the local cohomology of
graded modules over semigroup rings; in particular we prove that for a
semigroup $Q$ such that $Q^\sat$ is simplicial, and a finitely generated
module $M$ over $k[Q]$ which is graded by $\qgp$, the Bass numbers of
$H^i_I(M)$ are finite for any $Q$-graded ideal $I$ of $k[Q]$.
Conversely, if $Q^\sat$ is not simplicial, one can find a graded ideal
$I$ and graded $R$-module $M$ such that the local cohomology module
$H^i_I(M)$ has infinite-dimensional socle.  We introduce and exploit the
cominatorially defined {\em essential set} of a semigroup.
%\vskip 1ex
%\noindent
%{{\it AMS Classification:} ; }
\vspace{-.9em}
\end{abstract}
\maketitle

%%%%%%%%%%%%%%%%%%%%%%%%%%%%%%%%%%%%%%%%%%%%%%%%%%%%%%%%%%%%%%%%%%%%%%%%%
{}%%%%%%%%%%%%%%%%%%%%%%%%%%%%%%%%%%%%%%%%%%%%%%%%%%%%%%%%%%%%%%%%%%%%%%%
%%%%%%%%%%%%%%%%%%%%%%%%%%%%%%%%%%%%%%%%%%%%%%%%%%%%%%%%%%%%%%%%%%%%%%%%%

%%%%%%%%%%%%%%%%%%%%%%%%%%%%%%%%%%%%%%%%%%%%%%%%%%%%%%%%%%%%%%%%%%%%%%%%%
\section{Introduction}%%%%%%%%%%%%%%%%%%%%%%%%%%%%%%%%%%%%%%%%%%%%%%%%%%%

The local cohomology modules $H^i_I(M)$ for finitely generated
modules~$M$ over noetherian rings~$R$ have been studied for several
decades.  When $I$ is a maximal ideal of $R$ the local cohomology of~$M$
is fairly well-understood (and reasonably well-behaved), but for general
ideals~$I$ much less is known, and the behavior can be quite bad.  For
instance, Hartshorne \cite{HarCofinite} has shown that the Bass numbers
of~$H^i_I(M)$ need not be finite for general~$I$ and~$R$.

Recently, however, some progress has been made in special cases.
When~$R$ is a regular local ring containing a field,
Lyubeznik~\cite{Lyu1}, and Huneke and Sharp~\cite{HuSh} have shown
that~$H^i_I(M)$ has finite Bass numbers.  In the same spirit (albeit by
different techniques), Yanagawa~\cite{YanPoset} has shown that
if~$\omega$ is the canonical module of a simplicial and normal semigroup
ring and~$I$ is a monomial ideal, then~$H^i_I(\omega)$ has finite Bass
numbers.

Our approach is a substantial generalization of that found
in~\cite{YanPoset}.  We consider a noetherian ring~$R$ graded by a
semigroup~$Q$, and modules over~$R$ graded by~$\qgp$.  In this setting,
we introduce a functor called the {\em \v Cech hull}
(Section~\ref{sec:cech}), which allows us to recover the full local
cohomology of a finitely generated module~$M$ from the portion of the
local cohomology that lies in those graded degrees which are elements
of~$Q$ (Section~\ref{sec:localcoh}).  This piece of the local cohomology
is often easier to understand than the entire module; in particular it is
finitely generated when \hbox{$R = k[Q]$}.  Using this fact we prove
several finiteness results for the local cohomology of $\qgp$-graded
modules over semigroup rings (Section~\ref{sec:finiteness}).  In
particular we show (without Cohen-Macaulay hypotheses) that when~$Q^\sat$
is simplicial, $H^i_I(M)$ has finite Bass numbers.  We also prove the
converse, constructing a local cohomology module with
infinite-dimensional socle when~$Q^\sat$ is not simplicial
(Section~\ref{sec:infinite}).  The constructions are polyhedral in
nature, exploiting a new combinatorial structure, the {\em essential set}
of a semigroup (Section~\ref{sec:essential}).  Properties of the
essential set govern the associated primes, and to some extent the
module structure, of the local cohomology of the canonical module.

%\end{section}{Introduction}%%%%%%%%%%%%%%%%%%%%%%%%%%%%%%%%%%%%%%%%%%%%%
%%%%%%%%%%%%%%%%%%%%%%%%%%%%%%%%%%%%%%%%%%%%%%%%%%%%%%%%%%%%%%%%%%%%%%%%%
%%%%%%%%%%%%%%%%%%%%%%%%%%%%%%%%%%%%%%%%%%%%%%%%%%%%%%%%%%%%%%%%%%%%%%%%%
\section{The \v Cech Hull}\label{sec:cech}%%%%%%%%%%%%%%%%%%%%%%%%%%%%%%%

Let~$Q$ be a cancellative, commutative semigroup, and~$\qgp$ its
Grothendieck group, i.e.\ the group obtained from $Q$ by adjoining an
inverse for every element.  (The reader may safely assume for the
purposes of this paper that $Q$ is {\em affine}---that is, a finitely
generated submonoid of $\ZZ^d$ with $\qgp = \ZZ^d$; indeed, we will make
this assumption starting in Section~\ref{sec:semigroups}.  However, we
hope that the extra generality in this section and the next will be
useful for more general semigroup gradings, such as those arising in the
Cox homogeneous coordinate rings of toric geometry \cite{Cox}.)  We say a
ring $R$ is {\em $Q$-graded} and an $R$-module $M$ is $\qgp$-graded if we
are given direct sum decompositions
$$
  R = \bigoplus_{a \in Q} R_{a} \quad {\rm and} \quad M =
  \bigoplus_{\alpha \in \qgp} M_\alpha
$$
such that $R_{a}R_{b} \subseteq R_{a + b}$ and $R_{a}M_{\beta} \subseteq
M_{a + \beta}$.  The category of $\qgp$-graded modules is henceforth
denoted by $\ModR$.
 
A morphism $M \rightarrow N$ in $\ModR$ is a degree-preserving $R$-module
homomorphism; i.e.\ a map~$f$ of $R$-modules such that $f(M_\alpha)
\subseteq N_\alpha$.  We denote by $\hom_R(M,N)$ the $R_0$-module of such
morphisms and by $\hhom_R(M,N)$ the $\qgp$-graded $R$-module
$$
  \hhom_R(M,N) = \bigoplus_{\alpha \in \qgp} \hom_R(M,N(\alpha)) =
  \bigoplus_{\alpha \in \qgp} \hom_R(M(-\alpha),N)\,.
$$
Here, the $\qgp$-graded $R$-module $N(\alpha)$ is the {\em shift of $N$
by $\alpha$}, defined by $N(\alpha)_\beta = N_{\alpha + \beta}$.

For a subset $S \subseteq \qgp$ closed under the action of $Q$, we define
the {\em $S$-graded part} $M_S \subseteq M$ to be
$$
  M_S = \bigoplus_{\alpha \in S} M_{\alpha}\,.
$$
We frequently consider the case in which $S = \alpha+Q$ for some $\alpha
\in \qgp$.  In particular, taking~$M_Q$ yields the part of $M$ graded by
elements of $Q$.  Taking $S$-graded parts is functorial and exact, for
any $S$.

The theory we develop below revolves around the following question: To
what extent can we recover a module $M$ from its $Q$-graded part $M_Q$?
If $M$ is finitely generated, for example, then although we may not be
able to get $M$ from $M_Q$, we can shift by some $a \in Q$ to get $M(-a)
= M(-a)_Q$.  Therefore, the question is more meaningful for infinitely
generated modules, such as the local cohomology modules of a finitely
generated module.  We will find that these belong to a certain class of
modules that can be recovered from $Q$-graded parts of other modules by a
functor $\cech$ that we call the {\em \v Cech hull}, as defined by the
next result.

\begin{thm}
The functor $(-)_Q: \ModR \rightarrow \ModR$ taking $Q$-graded parts has
a right adjoint $\cech$; that is, there exists a functor $\cech$ and
natural isomorphisms
$$
  \hom_R(M_Q, N) = \hom_R(M, \cech N)
$$
for any $M,N \in \ModR$.
\end{thm}
\proof Given $N$, we explicitly construct $\cech N$ by defining
\begin{equation} \label{eqn:cech}
\begin{array}{r@{\ =\ }l}
    (\cech N)_{\alpha}	& \hom_R(R_{Q-\alpha}, N(\alpha)) \\
			& \hom_R(R_{Q-\alpha}(-\alpha), N))\,.
\end{array}
\end{equation}
The multiplication maps
$$
  R_{b} \otimes_{R_\0}(\cech N)_{\alpha} \rightarrow (\cech N)_{b+\alpha}
$$
are given by taking $r \otimes \phi$ to $(x \mapsto \phi(rx))$.  This is
well-defined since multiplication by $r \in R_{b}$ is a degree zero map
$R_{Q-\alpha-b}(-\alpha-b) \rightarrow R_{Q-\alpha}(-\alpha)$.

Note that if $a \in Q$ then $R_{Q-a} = R$, so
$$
  (\cech N)_{a} = \hom_R(R_{Q-a}(-a), N) = \hom_R(R(-a), N) = N_{a}\,,
$$
whence $(\cech N)_Q = N_Q$.  Therefore, given an element of
$\hom_R(M, \cech N)$, taking its $Q$-graded part gives an element of
$\hom_R(M_Q, N_Q)$.  This last module is $\hom_R(M_Q,N)$
(since degree zero maps from $M_Q$ to $N$ must land in $N_Q$),
so we have produced a natural map $\hom_R(M, \cech N) \to
\hom_R(M_Q,N)$.

Conversely, if $f \in \hom_R(M_Q, N)$ then for each $\alpha \in
\qgp$ we have a map
$$
  M_{\alpha} \rightarrow (\cech N)_{\alpha} = \hom_R(R_{Q-\alpha},
  N(\alpha)) \quad {\rm defined\ by} \quad x \mapsto (r \mapsto f(rx))\,.
$$
This is well-defined since if $r \in R_{Q-\alpha}$ and $x \in
M_{\alpha}$, then $rx \in M_Q$, so we can evaluate $f(rx)$.  We thus
obtain a well-defined element of $\hom_R(M, \cech N)$, whose $Q$-graded
part is just $f$.  This gives the natural inverse map for our bijection.
\endproof

\begin{remark} \rm
It is clear by looking at the graded pieces of $\cech$ that it is left
exact but not right exact, and that its derived functors are given in
terms of $\ext$ modules; that is $((R^i\cech)M)_{\alpha} =
\ext^i_R(R_{Q-\alpha}(-\alpha), M).$ Since $R_{Q-\alpha}(-\alpha)$ is
supported on $Q \subset \qgp$, both the {\v C}ech hull and its derived
functors depend only on the $Q$-graded part of $M$.
\end{remark}

\begin{remark} \rm \label{rk:poly}
The \v Cech hull was defined in \cite{Mil2} for polynomial rings.  In
this case $Q = \NN^d$, $R = k[Q]$, and $\hom_R(R_{Q -\alpha}(-\alpha), M)
= M_{\alpha^+},$ where $\alpha^+$ is obtained from $\alpha$ by zeroing
out the negative coordinates.  Thus our definition of $\cech$ agrees with
the one in \cite{Mil2}.  Note that the \v Cech hull is exact in this
case, since $R_{Q-\alpha}$ is free for all $\alpha$, so the $\ext$
modules that make up the graded pieces of $R^i\cech$ vanish.
\end{remark}

\begin{example} \rm \label{ex:cech}
Let $Q \subset \ZZ^2$ be the semigroup generated by $(0,2)$,$(1,1)$, and
$(2,0)$, so that $\qgp \subset \ZZ^2$ is a lattice of index 2.  Take
$R=k[Q]$, graded by $Q$.  If $\alpha = (x,y) \in \qgp$, then $x$ and $y$
have the same parity.  If $x$ and $y$ are even, or if $x$ and $y$ have
the same sign, then $R_{Q-\alpha}(-\alpha)$ is free.  On the other hand,
if $x$ is odd and negative while $y$ is odd and positive, then
$R_{Q-\alpha}(-\alpha)$ is generated in degrees $(1,y)$ and $(0,y+1)$.
Moreover, one has an exact sequence:
$$
  0 \to R_{Q-\alpha}(-\alpha-(1,1)) \too R(-1,-y) \oplus R(0,-y-1)
  \too R_{Q-\alpha}(-\alpha) \to 0.
$$
Splicing homological and graded shifts of this short exact sequence
together gives a free resolution $F_{\spot}$ of $R_{Q-\alpha}(-\alpha)$
such that $F_i = R(-1-i,-y-i) \oplus R(-i,-y-i-1)$.  Since $Q$ is
symmetric in $x$ and $y$, a similar result holds with $x$ and $y$
reversed.

We now have free resolutions of $R_{Q-\alpha}(-\alpha)$ for every $\alpha
\in \qgp$, and we can use them to compute the Cech hull and its derived
functors.  For instance, consider the module $k(-u,-v)$ consisting of a
single copy of the residue field $k$, supported in a degree $(u,v)$
satisfying $v > u > 1$.  Then for $\alpha = (x,y)$, we find that
$\hom_R(R_{Q-\alpha}(-\alpha),k(-u,-v))$ is only nontrivial if $\alpha =
(u,v)$.  Since Equation~(\ref{eqn:cech}) implies that
$$
  R^i\cech (k(-u,-v))_{\alpha} \cong \ext_R^i(R_{Q-\alpha}(-\alpha),
  k(-u,-v)),
$$
it follows that $\cech (k(-u,-v)) = k(-u,-v)$.  In contrast, $R^i\cech
(k(-u,-v))$ for $i > 1$ can only be nonzero in degrees $\alpha$ for which
$x$ and $y$ are odd and of differing sign, since $R_{Q-\alpha}(-\alpha)$
is free otherwise.

Suppose we have such an $\alpha$, and let $F_{\spot}$ be the
corresponding free resolution of $\alpha$ constructed above.  Then
$\hom_R(F_i, k(-u,-v))$ is nonzero if (and only if) one of the generators
of $F_i$ sits in degree $(u,v)$.  Referring to the expression for the
degrees of the $F_i$ in terms of $\alpha$ (and remembering that we
assumed $v > u > 1$), we find that $\hom_R(F_i,k(-u,-v))$ is nonzero if
and only if $i = u$ and $y = v-u-1$ or if $i = u-1$ and $y = v-u+1$.  In
other words, $R^i\cech (k(-u,-v))$ vanishes except when $i \in
\{u-1,u\}$.  Moreover, $R^{u-1}\cech (k(-u,-v))$ is supported in those
degrees $\alpha$ such that $x$ is odd and negative and $y = v-u+1$, while
$R^u\cech (k(-u,-v))$ is supported in those degrees $\alpha$ such that
$x$ is odd and negative and $y = v - u - 1$.

To summarize, when $v > u > 1$, we have:
\begin{itemize}
\item $\cech (k(-u,-v)) = k(-u,-v)$
\item $R^{u-1}\cech (k(-u,-v)) = k[\xx^{-(2,0)}](-1,-v + u - 1)$
\item $R^u\cech (k(-u,-v)) = k[\xx^{-(2,0)}](-1,-v + u + 1)$
\end{itemize}
and all other derived functors of $\cech$ vanish.  Here, $\xx^\alpha \in
k[\qgp]$ is the element corresponding to $\alpha \in \qgp$.
\end{example}

%\end{section}{The \v Cech hull}%%%%%%%%%%%%%%%%%%%%%%%%%%%%%%%%%%%%%%%%%
%%%%%%%%%%%%%%%%%%%%%%%%%%%%%%%%%%%%%%%%%%%%%%%%%%%%%%%%%%%%%%%%%%%%%%%%%
%%%%%%%%%%%%%%%%%%%%%%%%%%%%%%%%%%%%%%%%%%%%%%%%%%%%%%%%%%%%%%%%%%%%%%%%%
\section{Local Cohomology}\label{sec:localcoh}%%%%%%%%%%%%%%%%%%%%%%%%%%%

In this section we study the interaction of the \v Cech hull with the
functor $\Gamma_I$, which takes a module $M$ to the submodule annhilated
elementwise by some power of the ideal $I$.  If $I$ is graded, then
$\Gamma_I$ takes the category $\ModR$ of $\qgp$-graded modules to itself.
%$\ModR$.

\begin{prop} \label{prop:commute}
If $R$ is noetherian, then $\cech$ and $\Gamma_I$ commute: $\Gamma_I\cech
= \cech\Gamma_I$.
\end{prop}

\vbox{
\smallskip
{\noindent\it Proof.\/\,}

\vspace{-5ex}
$$
\begin{array}{r@{\ \: = \ \:}l}
(\Gamma_I\cech M)_\alpha
	& \bigcup_n\{x \in \hom_R(R_{Q-\alpha}(-\alpha), M) \mid I^nx = 0\}\\ 
	& \hom_R(R_{Q -\alpha}(-\alpha), \Gamma_IM)\\ 
	& (\cech\Gamma_I M)_\alpha.
\end{array}
$$

\vspace{-2ex}
\endproof}

\medskip
Henceforth we assume $R$ is noetherian.  In
Proposition~\ref{prop:spectral} we shall apply the spectral sequence of a
composite functor to the functors $\Gamma_I \cech$ and $\cech \Gamma_I$.
In order to do this we use the fact that both $\Gamma_I$ and $\cech$ take
injectives to injectives.  For $\Gamma_I$ this is standard; for $\cech$
this follows from Lemma~\ref{lemma:CJ}, whose extra precision is vital
for Section~\ref{sec:finiteness}.

\begin{lemma} \label{lemma:CJ}
Let $J$ be an indecomposable injective in $\ModR$.  Then $\cech J = 0$ if
$J_Q = 0$, and $\cech J = J$ otherwise.  In particular, $\cech$ takes
injectives to injectives.
\end{lemma}
\proof The first statement is clear, since $\cech$ depends only on the
$Q$-graded portion of a module.  The last statement follows from the
others, using the fact that every injective is a direct sum of
indecomposable injectives (this uses the noetherian hypothesis).

For the remaining statement, write $J = \EE(R/\pp)(\alpha)$ for some
prime $\pp$ of $R$ and $\alpha \in \qgp$.  Then $(R/\pp)(\alpha)$ is an
essential submodule of $J$, so since $J_Q$ is nonzero, 
$(R/\pp)(\alpha)_Q$ is nonzero.  Now $(R/\pp)(\alpha)_Q$ is an essential 
submodule of $(R/\pp)(\alpha)$ (since $R/\pp$ is an integral domain), so it 
is an essential submodule of $J$.  Thus in particular $J_Q$ is an essential
submodule of $J$.

The inclusion $J_Q \rightarrow J$ induces a map $\phi: J \rightarrow
\cech J$, by the adjointness property of $\cech$.  Moreover, $\phi$ is
injective, since it restricts to the identity on the essential submodule
$J_Q$.  Thus $J$ is a direct summand of $\cech J$, since $J$ is
injective.  We claim that $\cech J$ is an essential extension of $J$,
from which the result follows immediately.
 
Let $x \in (\cech J)_{\alpha} = \hom_R(R_{Q-\alpha}(-\alpha), J_Q)$
be a nonzero homogeneous element.  Then there is an $r \in R$ such that
$x(r)$ is nonzero, and then $rx = x(r)$ is a nonzero element of
$J_Q$.  Thus $\cech J$ is an essential extension of $J_Q$ and
hence of $J$ as well.
\endproof

\begin{prop} \label{prop:spectral}
Let $M$ be a graded $R$-module.  There are spectral sequences $E(M)$ and
$F(M)$ described by
$$
\begin{array}{r@{\ \: = \ \:}l}
  E_2^{p,q}(M) & R^p\cech H^q_I(M) \Rightarrow R^{p+q}(\Gamma_I \cech)M\\
  F_2^{p,q}(M) & H^p_I(R^q\cech M) \Rightarrow R^{p+q}(\Gamma_I \cech)M.
\end{array}
$$
% The filtrations of $R^{p+q}(\Gamma_I \cech)M$ to which they converge
% are not necessarily equal.
\end{prop}
\proof These are spectral sequences for the composite functors $\cech
\Gamma_I$ and $\Gamma_I \cech$ \cite{Wei}, using
Proposition~\ref{prop:commute} along with Lemma~\ref{lemma:CJ} and the
fact that $\Gamma_I$ takes injectives to injectives.  \endproof

\begin{example} \rm \label{ex:poly}
We return to the case of Remark~\ref{rk:poly}; that is, $R=k[\NN^d]$.
Here $\cech$ is exact, and so the spectral sequences $E(M)$ and $F(M)$
both collapse.  The proposition simply says that $H_I^i(\cech M) = \cech
H_I^i(M)$.  The right derived functors of the {\v C}ech hull measure the
degree to which this equality fails in other rings.
\end{example}

\medskip

Ultimately, the goal of this section is Theorem~\ref{thm:spectral}, which
describes how the local cohomology modules of $M$ can be reconstructed
from a finite collection of submodules thereof, using the \v Cech hull
and its derived functors.  After a choosing a certain $\alpha \in \qgp$,
this is accomplished by a spectral sequence $E(M(-\alpha))$ that only
depends on the (finitely generated) $Q$-graded parts of the local
cohomology modules of $M(-\alpha)$, since $R^p\cech H^q_I(-)$ depends
only on the $Q$-graded part of $H_I^q(-)$.

%Of course, we need to know that $R^{p+q}(\Gamma_I\cech)M$ is a local
For this approach to work, we of course need $R^{p+q}(\Gamma_I\cech)M$ to
be a local cohomology module.  To this end, we use the spectral sequence
$F(M)$.  Although the filtration that arises from $F(M)$ is generally
nontrivial, we avoid this nuisance by replacing $M$ with a suitable
$\qgp$-graded shift, forcing $F$ to collapse in low cohomological degree.
We find this suitable shift in Corollary~\ref{cor:shift} using
Proposition~\ref{prop:pinnacle}, which is interesting in its own right.

\begin{prop} \label{prop:pinnacle}
Let $J^\spot$ be a minimal injective resolution of a finitely generated
module $M \in \ModR$.  Let $\pp$ be a homogeneous prime of $R$, let
$\mm$ be a homogeneous maximal ideal containing $\pp$, and let 
$c=\dim(R/\pp) - \dim(R/\mm)$.  If every indecomposable
summand of $\Gamma_\mm J^{i+c}$ has nonzero $Q$-graded part, then every
indecomposable summand of $J^i$ isomorphic to a shift of $\EE(R/\pp)$
has nonzero $Q$-graded part.
\end{prop}
\proof Inverting all homogeneous elements outside $\mm$ fixes all shifts
of $\EE(R/\pp)$ as well as $\Gamma_\mm J^i$, so we assume henceforth that
$\mm$ is the unique maximal homogeneous ideal of $R$.

We begin with the case $c = 1$.  Using $(-)_{(\pp)}$ to denote the
localization by all homogeneous elements outside of $\pp$, it is a
standard fact (in \cite[p.~101]{BH}, for instance) that
$(\eext^i_R(R/\pp,M))_{(\pp)} = (\hhom_R(R/\pp,J^i))_{(\pp)} \subset
(\Gamma_\pp J^i)_{(\pp)}$ is an essential extension.  Therefore, we need
only show that every indecomposable submodule of the free
$(R/\pp)_{(\pp)}$-module $(\eext^i_R(R/\pp,M))_{(\pp)}$ has nonzero
$Q$-graded part.

Choose a homogeneous element $x \in \mm \setminus \pp$, and define $L$ by
the exact sequence
$$
  0 \to R/\pp \stackrel{x}{\too} R/\pp \too L \to 0,
$$
where $\stackrel{x}{\to}$ is multiplication by $x$.  The long exact
sequence for $\eext^\spot_R(-,M)$ provides a right exact sequence
$$
  \eext^i_R(R/\pp,M) \stackrel{x}{\too} \eext^i_R(R/\pp,M) \too
  L(i,x,M) \to 0
$$
for the appropriate submodule $L(i,x,M) \subseteq \eext^{i+1}_R(L,M)$.
Tensoring with $R/\mm$ yields an isomorphism $R/\mm \otimes
\eext^i_R(R/\pp,M) \cong R/\mm \otimes L(i,x,M)$ since multiplication by
$x$ becomes the zero map.  Nakayama's lemma implies that
$\eext^i_R(R/\pp,M)$, and hence $(\eext^i_R(R/\pp,M))_{(\pp)}$, is generated
by elements in degrees $\gamma \in \qgp$ such that $L(i,x,M)_\gamma
\subseteq \eext^{i+1}_R(L,M)_\gamma \neq 0$.

Now $\eext^{i+1}_R(L,M)$ is the $(i+1)^{\rm st}$ cohomology of the complex
$\hhom_R(L,J^\spot) \subset \Gamma_\mm J^\spot$, whose socle subcomplex
$\hhom_R(R/\mm,J^\spot) \subseteq \hhom_R(L,J^\spot)$ (the inclusion being
induced by the surjection $L \twoheadrightarrow R/\mm$) is equal to
$\eext^\spot_R(R/\mm,M)$.  The hypothesis on $(\Gamma_\mm J^i)_Q$ in the
Proposition implies that $\eext^{i+1}_R(R/\mm,M)$ must in fact equal its
$Q$-graded part.  Given any degree $\gamma \in \qgp$ for which
$\eext^{i+1}_R(L,M)_\gamma \neq 0$, we therefore can find a homogeneous element
$r \in R$ such that $\gamma + \deg(r) \in Q$.  Note that $\pp$
annihilates $L$ and hence also $\hhom_R(L,J^\spot)$, so we can always find
our element $r$ outside of $\pp$.  When we invert $r$ to form the
localization $(\eext^i_R(R/\pp,M))_{(\pp)}$, any generator $y$ in degree
$\gamma$ can be replaced by the generator $ry$ whose degree is in $Q$.
This concludes the case where $c = 1$.

The general case proceeds by induction on $c$, replacing $R$ with its
homogeneous localization at a prime containing $\pp$ and having dimension
$\dim(R/\pp) - 1$.~\endproof

\medskip
Recall that the {\em Bass number}
%$\mu_i(\pp,M)$}
of a module $M$ at the prime $\pp$ in cohomological degree $i$ is the
number of indecomposable summands isomorphic to a shift of $\EE(R/\pp)$
appearing at the $i^{\rm th}$ stage in any minimal injective resolution
of $M$.  These numbers are always finite if $M$ is finitely generated,
but may in general be infinite.
% This is also equal to the vector space dimension of the homogeneous
% localization $\eext^i_R(R/\pp,M)_{(\pp)}$ as an
% $(R/\pp)_{(\pp)}$-module, or simply the rank of $\eext^i_R(R/\pp,M)$ as
% an $R/\pp$-module.
\begin{cor} \label{cor:shift}
Suppose $R$ has a unique maximal homogeneous ideal $\mm$.
% (so $R_\0$ is local).
Let $M \in \ModR$ be a finitely generated $R$-module,
and $n$ be a positive integer.  Then there exists
$\alpha \in \qgp$ such that for all $\beta \in \alpha + Q$,
\begin{enumerate}
\item
$\cech (M(-\beta)) = M(-\beta)$, and

\item
$R^j \cech (M(-\beta)) = 0$ if $1 \leq j < n.$ 
\end{enumerate}
\end{cor}
\proof Since $M$ is finitely generated, the Bass numbers of $M$ at $\mm$
are finite.  Thus $\Gamma_\mm J^i$ is a finite direct sum of
indecomposables for each $i \leq n + \dim(R)$.  These can be moved to
have nonzero $Q$-graded part by some shift $(-\alpha)$.
Lemma~\ref{lemma:CJ} and Proposition~\ref{prop:pinnacle} together then
imply that $\cech$ fixes $J^\spot(-\beta)$ in cohomological degree $n$
and less for all $\beta \in \alpha + Q$.~\endproof

\begin{remark} \label{remark:shift} \rm
If $R$ is a ring with only finitely many homogeneous primes (e.g.\ a
semigroup ring), then the conclusion of Corollary~\ref{cor:shift} holds
for any $M$ with finite Bass numbers, as then $M$ has an injective
resolution with finitely many summands in each cohomological~degree.
\end{remark}

\begin{example} \rm \label{ex:can}
Let $Q \subset \ZZ^d$ be a finitely generated semigroup, and let $R =
k[Q]$, graded by $Q$.  Ishida \cite{Ish} constructed a dualizing complex
for $R$, in which each indecomposable injective appears without shift.
When $R$ is Cohen-Macaulay, this is an injective resolution of the
canonical module $\omega_R$ that is fixed by the \v Cech hull.  Hence
Corollary~\ref{cor:shift} holds for $\omega_R$ with $\alpha = 0$.  It
follows that for $q > 0$ we have
$$
  F_2^{p,q}(\omega_R) = H^p_I(R^q\cech \omega_R) = 0,
$$
so $F(\omega_R)$ converges to $H^{p+q}_I(\omega_R)$.  Thus $E(\omega_R)$
likewise converges to a filtration of $H^{p+q}_I(\omega_R)$.  We will use
this fact in the next section (Proposition~\ref{prop:canonical}) to
compute $H^{p+q}_I(\omega_R)$ for $I$ prime.
\end{example}

\begin{example} \rm
The phenomenon predicted by Corollary~\ref{cor:shift} is clearly
illustrated in Example~\ref{ex:cech}: as $u$ and $v$ increase, the
derived functors of $\cech (k(-u,-v))$ in positive cohomological degrees
$< u-1$ vanish.
\end{example}

\begin{thm} \label{thm:spectral}
Suppose $R$ has a unique homogeneous maximal ideal.
Let $M$ be a finitely generated $R$-module, and $n$ be a positive
integer.  Then there exists $\alpha \in \qgp$ such that for all $\beta
\in \alpha + Q$, the spectral sequence
$$
  E_2^{p,q}(M(- \beta)) = R^p\cech H^q_I(M(- \beta))
  \Rightarrow H^{p+q}_I(M)(- \beta)
$$
converges to a local cohomology module for $p + q < n$.  \comment{this
result is likely true without any locality hypothesis on the nature of
the zeroth graded piece.  all one really needs is that there are finitely
many ``worst'' closed points in the spectrum of the zeroth graded piece.
one ought to be able to stratify the spectrum of the zeroth graded piece
into smooth pieces in such a way that the Bass degrees are bounded on
each piece, so that one forceful shift gives everything a nonzero
Q-graded part.  in any case, the only case I can think of where this
might be useful is when the zeroth graded piece is the ring of integers;
in almost all applications, the zeroth graded piece is a field, or at the
very least, local.}
\end{thm}
\proof Choose $\alpha$ as in Corollary~\ref{cor:shift}.  Then for all
$\beta \in \alpha + Q$,
$$
  F_2^{p,q}(M(- \beta)) = H_I^p(R^q\cech M(- \beta))
  = \left\{
	\begin{array}{cc}
	0  & \text{if }q>0\\
	H_I^p(M)(- \beta) & \text{if } q=0
	\end{array}\right. .
$$
Hence if $p + q < n$, $R^{p+q}(\Gamma_I\cech)(M(- \beta)) =
H_I^{p+q}(M(- \beta))$.  Since $E(M(- \beta))$ converges
to the former by Proposition~\ref{prop:spectral}, the result follows.
\endproof

Injective resolutions are rarely finite, so no matter which $\alpha$ is
chosen in Theorem~\ref{thm:spectral}, $F(M)$ really can converge to
something other than $H^{p+q}_I(M)$ in large cohomological degrees.
For example, if we take $Q$, $R$, and $k(-u,-v)$ as in Example~\ref{ex:cech},
then $H^i_I(k(-u,-v))$ vanishes for $i \geq 1$ and any $I$.  On the other
hand, we have $F_2^{p,q}(k(-u,-v)) = 0$ for $p \geq 1$ and $F_2^{0,q}$ nonzero
for $q \in \{u-1,u\}$, so the nonvanishing derived functors of $\cech$
cause $F(k(-u,-v))$ to fail to converge to local cohomology in these degrees. 

However, since $H^{p+q}_I(M)$ vanishes in sufficiently high cohomological
degrees, choosing $n$ large in the theorem does show how the collection
of $Q$-graded parts $H^j_I(M)(-\beta)_Q$ for all $j$ determine the entire
local cohomology modules.  As we shall see in
Section~\ref{sec:semigroups}, the $Q$-graded portion of a local
cohomology module is often much easier to understand than the local
cohomology module itself.

%\end{section}{Local Cohomology}%%%%%%%%%%%%%%%%%%%%%%%%%%%%%%%%%%%%%%%%%
%%%%%%%%%%%%%%%%%%%%%%%%%%%%%%%%%%%%%%%%%%%%%%%%%%%%%%%%%%%%%%%%%%%%%%%%%
%%%%%%%%%%%%%%%%%%%%%%%%%%%%%%%%%%%%%%%%%%%%%%%%%%%%%%%%%%%%%%%%%%%%%%%%%
\section{Semigroup Rings}\label{sec:semigroups}%%%%%%%%%%%%%%%%%%%%%%%%%%

One of the ways of understanding local cohomology $H^\spot_I(-)$ in terms
of finitely generated modules is by taking limits (over $m$) of modules
$\ext^\spot_R(R/I^m, -)$.  Unfortunately, these limits are frequently
quite badly behaved (see \cite{EMS}, for example)\comment{other
references?}.  Here, we bypass them entirely, in the case where $R =
k[Q]$ is an affine semigroup algebra over a field $k$, by constructing
the graded pieces of $H^\spot_I(M)$ in terms of the derived functors
$\eext^\spot_R(R/I^m,M)$ of $\hhom_R(R/I^m,M)$ for a {\em single fixed}
$m$, using Theorem~\ref{thm:spectral}.  In order for this to work, we
need to know what the $Q$-graded part of local cohomology looks like.

In Sections~\ref{sec:semigroups}--%
% ,~\ref{sec:finiteness},~\ref{sec:essential}, and~
\ref{sec:infinite}, we set $R=k[Q]$, an affine semigroup algebra graded
by $Q \subseteq \ZZ^d$, which is not assumed normal.  Such a ring
satisfies the hypotheses of Corollary~\ref{cor:shift}, so that all of the
machinery of the previous sections applies.  For $k[Q]$ we also have a
simpler expression for the \v Cech hull.  In what follows, $Q$ is viewed
as contained in $k[Q]$ via $a \mapsto \xx^a$.
\begin{prop} \label{prop:cechsemigroup}
When $R=k[Q]$ and $M \in \ModR$, we have $(\cech M)_{\beta} \cong
\hom_R(R_{Q+\beta},M)$.  If $a \in Q$, we have a commutative diagram:
$$
\begin{array}{c@{\ }c@{\ }c}
	(\cech M)_{\beta}
&	\stackrel{\cdot \textstyle \xx^a}\longrightarrow
&	(\cech M)_{\beta+a}
\\
	\downarrow
&
&	\downarrow
\\
\hom_R(R_{Q+\beta},M)
&	\longrightarrow
&	\hom_R(R_{Q+a+\beta},M),
\end{array}
$$
where the vertical arrows are isomorphisms and the bottom arrow
is induced by the inclusion of $R_{Q+a+\beta}$ in $R_{Q+\beta}$.
\end{prop} 
\proof By definition, $(\cech M)_{\beta} = 
\hom_R(R_{Q-\beta}(-\beta), M).$  Multiplication by $\xx^\beta$ 
induces an injection $R_{Q-\beta}(-\beta) \rightarrow R_{Q+\beta}$;
this is an isomorphism since both of these modules are supported in the
same degrees.   
 
Moreover, one has the following commutative diagram:
$$
\begin{array}{c@{\ }c@{\ }c}
R_{Q-\beta-a}(-\beta-a)
&	\stackrel{\cdot \textstyle \xx^a}{\longrightarrow}
&	R_{Q-\beta}(-\beta)
\\
\downarrow
&
&	\downarrow
\\
R_{Q+\beta+a}
&	\longrightarrow
&	R_{Q+\beta}
\end{array}
$$
from which the rest of the proposition follows immediately. 
\endproof

\begin{lemma} \label{lemma:ext}
Let $J = \EE(R/\pp)(-\alpha)$ be an indecomposable injective, and $I$ an
ideal of~$R$.  There exists $n \in \NN$ such that for all $m > n$,
$(\Gamma_I J)_Q = \hhom_R(R/I^m, J)_Q.$
\end{lemma}

\proof Suppose $I$ is not contained in $\pp$.  Then some element of $I$
acts as a unit on $R/\pp$, so $\Gamma_I J = \hhom_R(R/I^m, J) = 0$ and
the result is trivial.  Thus it suffices to show this result for $I$
contained in $\pp$.  In this case $\Gamma_I J = J$, so it suffices to show
that $\hhom_R(R/I^m, J)_Q = J_Q$; i.e.\ that every element of $J_Q$ is
killed by $I^m$.

Let $\tau$ be a linear functional that takes nonnegative values on $Q$,
such that if $b \in Q$, then $\tau(b) > 0 \Leftrightarrow b \in \pp$.
Then $\EE(R/\pp)$ is supported in those degrees $\beta$ such that
$\tau(\beta) \leq 0$.  Thus $J$ is supported in those degrees $\beta$
such that $\tau(\beta) \leq \tau(\alpha) =: n$.

Suppose $m > n$.  Let $y \in J_Q$ and $x \in I^m$ be nonzero homogeneous
elements of degrees $b$ and $c$, respectively.  Then $\tau(b) \geq 0$
because $b \in Q$ and $\tau(c) \geq m$ because $x \in I^m \subseteq
\pp^m$.  Thus $xy$ lies in degree $b + c$ and $\tau(b+c) > n$, so $xy =
0$, as required.  \endproof

\medskip
Now we can apply the Lemma to describe $Q$-graded parts of local
cohomology.
\begin{prop} \label{prop:ext}
Let $M \in \ModR$ be finitely generated, and $I$ be a graded ideal of
$R$.  Fix a nonnegative integer $i$.  Then there exists $m_0 \in \NN$
such that for any $m \geq m_0$,
$$
  H_I^i(M)_Q \cong \eext^i_R(R/I^m,M)_Q.
$$
\end{prop}

\proof Let $J^\spot$ be an injective resolution for $M$, and choose $m_0$
sufficiently large that $(\Gamma_I N)_Q = \hhom_R(R/I^m, N)_Q$ agree for
every $m \geq m_0$ and every inedecomposable injective summand $N$
appearing in cohomological degree $i$ or lower in $J^\spot$.  Then the
first $i$ right derived functors of $(\Gamma_I -)_Q$ and $\hhom_R(R/I^m,
-)_Q$ agree on $M$; since $(-)_Q$ is exact this means $H_I^i(M)_Q \cong
\eext^i_R(R/I^m, -)_Q$.  \endproof

\begin{example} \rm \label{ex:can'}
If $Q$ is saturated and $M = \omega_R$ then the power of $I$ in
Proposition~\ref{prop:ext} can be set equal to $1$; i.e.
$$
  \eext^p_R(R/I,\omega_R)_Q \cong H^p_I(\omega_R)_Q,
$$
since the indecomposable summands of the injective resolution of
$\omega_R$ are unshifted~$\EE(R/\pp)$'s.% in degree zero.
\end{example}

\begin{cor} \label{cor:qfinite}
If $M \in \ModR$ is finitely generated, $H^i_I(M)_Q$ is 
finitely generated.
\end{cor}

\begin{cor} \label{cor:cechfinite} 
Let $M \in \ModR$ be finitely generated.  Then $H^i_I(M)$ has a finitely
generated essential submodule if and only if for some $\beta \in \qgp$
the natural map $H^i_I(M)(-\beta) \rightarrow \cech(H^i_I(M)(-\beta))$ is
an injection.
\end{cor}
\proof Suppose $H^i_I(M)$ has a finitely generated essential submodule
$N$.  Then we can shift $N$ so that all of its generators have degrees in
$Q$; i.e.\ there exists $\beta$ such that $N(-\beta) \subset
(H^i_I(M)(-\beta))_Q$.  Since the map $H^i_I(M)(-\beta) \rightarrow \cech
H^i_I(M)(-\beta)$ is injective on its $Q$-graded part, it is injective on
$N(-\beta)$; since $N(-\beta)$ is essential the map is injective
everywhere.

Conversely, $H^i_I(M)(-\beta))_Q$ is an essential submodule of
$\cech(H^i_I(M)(-\beta))$ and hence of $H^i_I(M)(-\beta)$.  By
Corollary~\ref{cor:qfinite} it is finitely generated.  \endproof

\medskip
The upshot of the above is that since the spectral sequence $E(M)$
depends only on the $Q$-graded parts of the local cohomology modules which
appear in it, we can just replace these local cohomology modules with
the corresponding $\eext$-modules.

\begin{thm} \label{thm:extspectral}
Let $M$ be a finitely generated module over $R = k[Q]$, and $n$ be a
positive integer.  Then there exists $\alpha \in \qgp$ such that for all
$\beta \in \alpha + Q$, there exists $m \in \ZZ$ making the spectral
sequence
$$
  E_2^{p,q}(M(- \beta)) = R^p\cech \eext_R^q(R/I^m, M(-\beta))
  \Rightarrow H^{p+q}_I(M)(- \beta)
$$
converge to a local cohomology module for $p + q < n$.  Taking degree
$\gamma$ parts for any $\gamma \in Q$ yields a spectral sequence of
iterated Ext modules,
$$
  E_2^{p,q}(M(-\beta))_\gamma = \ext_R^p(R_{Q+\gamma}(\beta),
  \eext_R^q(R/I^m, M)) \Rightarrow H^{p+q}_I(M)_{\gamma-\beta}\,.
$$
\end{thm}
\proof This is immediate from Theorem~\ref{thm:spectral} and
Proposition~\ref{prop:ext}.  \endproof

\begin{example}\label{ex:poly2} \rm
Returning to the setting of Example~\ref{ex:poly}, we find using this
theorem that if $M$ is a $\ZZ^d$-graded module over a polynomial ring in
$d$ variables and $I$ is a monomial ideal, then there exist $m \in \ZZ$
and $\beta \in \ZZ^d$ such that $H_I^i(M) = \cech \eext^i_R(R/I^m,
M(-\beta))$.  This generalizes a result proved independently by
Musta\c{t}\v{a} \cite{Mus1} and Terai \cite{TerLocalCoh}.
\end{example}

For saturated $Q$, Theorem~\ref{thm:extspectral} takes an especially nice
form for canonical modules.  Recall that a {\em face} of $Q$ is the set
of degrees of elements outside a prime ideal of $R$.

\begin{prop} \label{prop:canonical}
Suppose $R$ is normal and of dimension~$d$.  Let $\pp$ be a prime of $R$,
corresponding to an $n$-dimensional face of $Q$.  Then
$H^{d-i}_\pp(\omega_R) \cong R^{n - i}\cech(\omega_{R/\pp})$.
\end{prop}
\proof $R^q\cech\eext^p_R(R/\pp,\omega_R) \Rightarrow
H^{p+q}_\pp(\omega_R)$ by Example~\ref{ex:can'} and
Theorem~\ref{thm:extspectral}.  Since $R/\pp$ is a dimension~$n$
Cohen-Macaulay quotient of the Cohen-Macaulay ring $R$ of dimension~$d$,
the module $\eext^p_R(R/\pp,\omega_R)$ is nonzero only when $p = d-n$, in
which case it is $\omega_{R/\pp}$.
%(proof: $H^{n-p}_\mm(R/\pp)$ only vanishes for $p = n
%- \dim F$ because $R/\pp$ is Cohen-Macaulay of dimension $\dim F$, and
%the Matlis duals of these modules are the Ext modules in question).
Thus the spectral sequence degenerates, and $R^q\cech(\omega_{R/\pp})
\cong H_\pp^{q + d - n}(\omega_R)$.  \endproof

\begin{example} \label{ex:har} \rm
Let $Q$ be the semigroup on four generators $\{x,y,u,v\}$ and one relation
$x+u = y+v$, and $R=k[Q]$.  In \cite{HarCofinite}, Hartshorne shows that
for the ideal $I=(\xx^u,\xx^v)$, the local cohomology module 
$H^1_I(\omega_R)$ has a
finitely generated essential submodule while $H^2_I(\omega_R)$ has an
infinite dimensional socle, supported in degrees $n(x-v)$ for $n > 0$.  
This is consistent with
Proposition~\ref{prop:canonical}, which says that infinite-dimensional
socles must arise from a nonvanishing higher derived functor of the \v
Cech hull.  See Section~\ref{sec:infinite} for a combinatorial
explanation of why this bad behavior occurs, in the context of its
generalization to arbitrary affine semigroup rings.\comment{perhaps make
a comment about the cohomological dimension cd(I) in the gorenstein case,
cf. the comment below.  the discrepancy from cd(I) equalling the
cohomological degree of the last nonvanishing $\ext^i(R/I,R)$ is given by
the right derived functors of \v cech hull}
\end{example}
\comment{respond to \cite[Theorem~5.12]{YanPoset} by removing the
simplicial hypothesis.  give \cite[Theorem~5.9]{YanPoset} as an example
here, noting that we know much more now---specifically how to recover the
local cohomology modules from the $\ext$ modules, even when some of the
latter are mysteriously zero}

%\end{section}{Semigroup Rings}%%%%%%%%%%%%%%%%%%%%%%%%%%%%%%%%%%%%%%%%%%
%%%%%%%%%%%%%%%%%%%%%%%%%%%%%%%%%%%%%%%%%%%%%%%%%%%%%%%%%%%%%%%%%%%%%%%%%
%%%%%%%%%%%%%%%%%%%%%%%%%%%%%%%%%%%%%%%%%%%%%%%%%%%%%%%%%%%%%%%%%%%%%%%%%
\section{Finiteness for simplicial semigroups}\label{sec:finiteness}%%%%%

Letting $Q \subseteq \ZZ^d$ be affine as in the previous section, the
above machinery allows us to make strong statements about local
cohomology over $R=k[Q]$.  Indeed, the fact that local cohomology modules
``come from'' the derived functors of the \v Cech hull forces certain
structure on them.  This structure is codified in the notion of a
\emph{straight} module, a common generalization of notions due to Miller
\cite{Mil2} (who defined {\em $\mathbf a$-determined modules} over a
polynomial ring for $\mathbf{a} \in \NN^n$) and Yanagawa \cite{YanPoset,
Yan2} (who defined straightness for a restrictive class of modules over
semigroup rings).  One obtains from Theorem~\ref{thm:extspectral} that
local cohomology modules over semigroup rings are straight when shifted
appropriately.  Over a simplicial (and not necessarily normal) semigroup
ring this forces them to have finite Bass numbers
(Theorem~\ref{thm:finite}).  The key to all this is the following
definition.

\begin{defn} \sl \label{defn:straight} 
Let $Q \subseteq \ZZ^d$ be affine and $R= k[Q]$.  The category $\SS$ of
\emph{straight} $R$-modules is the smallest subcategory of the
$\qgp$-graded modules $\ModR$ such that (1) all indecomposable injectives
$J$ satisfying $J_Q \neq 0$ are in $\SS$; (2) finite direct sums of
modules in $\SS$ are in $\SS$; and (3) if $\phi$ is a homomorphism of
straight modules, then $\ker(\phi)$ and $\coker(\phi)$ are straight.
\end{defn}

%\begin{remark} \rm
%This notion first appeared in \cite{YanPoset}.  Yanagawa's definition was
%slightly different; as he was dealing with squarefree modules, he only
%considered those indecomposable injectives whose $Q$-graded parts were
%squarefree to be straight.  For our purposes this is too restrictive,
%hence the new definition.
%\end{remark}

\begin{prop} \label{prop:cechstraight}
Let $M \in \ModR$ be finitely generated.  Then $R^i\cech M$ is straight. 
\end{prop}
\proof By definition, a finite direct sum of straight modules is
straight.  Thus if $C^\spot$ is a complex of finite direct sums of
indecomposable injectives, and each indecomposable injective in $C^\spot$
has nontrivial $Q$-graded part, then $C^\spot$ is a complex of straight
modules, and the cohomology of $C^\spot$ is straight.  In particular, if
$C^\spot = \cech J^\spot$ is the \v Cech hull of an injective resolution
$J^\spot$ of $M$, then $\cech J^\spot$ is a complex of straight modules
by Lemma~\ref{lemma:CJ}.  Thus its cohomology is straight, as required.
\endproof

\begin{cor} \label{cor:locstraight}
Let $M \in \ModR$ be finitely generated.  Then there exists an element
$\,a \in Q$ such that $H^i_I(M)(-a)$ is straight.
\end{cor}
\proof By Theorem~\ref{thm:extspectral} and
Proposition~\ref{prop:cechstraight}, we have $a \in Q$ such that
$H^i_I(M)(-a)$ is the limit of a spectral sequence of straight modules.
This spectral sequence yields a {\em finite} filtration of $H^i_I(M)(-a)$
whose associated graded modules are straight.  \endproof

\medskip
Straight modules have a number of useful properties.  In particular, the
fact that they can be ``built out of'' indecomposable injectives with
nontrivial $Q$-graded part by taking kernels, cokernels, and finite
direct sums forces many of their graded pieces to be isomorphic to
each other.

\begin{prop} \label{prop:straight}
Let $M$ be a straight module over $R = k[Q]$.  Then:
\begin{enumerate}
\item
$M_Q$ is finitely generated.

\item
$M(-a)$ is straight for all $a \in Q$.

\item
Multiplication by $\xx^a$ is an isomorphism $M_\beta \rightarrow M_{a +
\beta}$ whenever $\beta \in \qgp$ and $\,a \in Q$ satisfy $(\beta+Q) \cap
Q = (a+\beta+Q) \cap Q$.
\end{enumerate}
\end{prop}
\proof If the above three properties hold for $M$ and $N$, then they also
hold for $M \oplus N$, as well as $\ker(\phi)$ and $\coker(\phi)$ for any
$\phi: M \to N$.  Thus it suffices to check that if $J$ is an
indecomposable injective, and $J_Q$ is nonzero, then $J$ has the above
properties.  $J_Q$ is clearly finitely generated, and if $J_Q$ is
nonzero, so is $J(-a)_Q$, so the first two properties are clear.

For the third property, note that $\cech J = J$ by Lemma~\ref{lemma:CJ}.
Thus in particular,
$$
\begin{array}{cl@{\ \: =\ \:}l@{\ \: =\ \:}l}
&	\:J_{\beta} 
&	(\cech J)_\beta
&	\hom_R(R_{Q+\beta}, J)
\\[5pt]
	\hbox{and}
&	J_{a + \beta}
&	(\cech J)_{a + \beta}
&	\hom_R(R_{Q+a+\beta}, J).
\end{array}
$$ 
% But if $\tau(\beta)^+ = \tau(\beta+\alpha)^+,$
The hypothesis in part~3 says that $R_{Q+\beta} = R_{Q+a+\beta}$, whence
multiplication by $\xx^a$ is an isomorphism $J_\beta \rightarrow
J_{a+\beta}$, as required.  \endproof

\medskip
The third property of Proposition~\ref{prop:straight} motivates the
following definition.

\begin{defn} \label{defn:essential}
An {\em essential point for $Q$} is an element $\varepsilon \in \qgp$
such that if $a \in Q$ and $(\varepsilon+Q) \cap Q = (a+\varepsilon+Q)
\cap Q$, then $\:a\:$ is a unit.  The {\em essential set} $\CE$ is the
$Q$-set generated by the essential points; i.e.\ $\CE$ is the union
$\bigcup (Q+\varepsilon)$ over essential points $\varepsilon$.
% (i.e. x is maximal in the Q-ordering among those points y such that y+Q
% \cap Q = x+Q \cap Q).
\end{defn}

We will use combinatorial properties of the essential set $\CE$ to prove
the main result of~this section, Theorem~\ref{thm:finite}, which is
algebraic.  However, we postpone the proofs of the combina\-torics to
Section~\ref{sec:essential}, where we investigate the structure of the
essential set in more detail.
\begin{prop}\label{prop:essential}
If $M$ is straight, then $M_\CE$ is an essential submodule of $M$.
\end{prop}
\vspace{-.75pt}
\proof Suppose $0 \neq x \in M_\alpha$, and choose an essential point
$\varepsilon$ with~\hbox{$\varepsilon-\alpha \in Q$} and $(\varepsilon +
Q) \cap Q = (\alpha + Q) \cap Q$, using Lemma~\ref{lemma:exists}.
Proposition~\ref{prop:straight} shows that multiplication by
$\xx^{\varepsilon-\alpha}$ is an isomorphism $M_\alpha \to
M_\varepsilon$.  Thus any submodule of $M$ containing $x$ contains
$\xx^{\varepsilon-\alpha}x \in M_\CE$.~\endproof
 
\medskip
Observe that the results in this section so far have used no extra
hypotheses on the affine semigroup $Q$.  Since our goal involves
simplicial semigroups, this will change starting now.
\begin{prop}\label{prop:strfinite}
Let $Q$ be an affine semigroup whose saturation is simplicial modulo
units, and $M$ a straight module over $k[Q]$.  Then the Bass numbers of
$M$ are finite.
\end{prop}
\proof Taking $a$ as in Proposition~\ref{prop:simplicial}, we find that
$M_\CE \subset M_{Q-a}$, so $M_\CE(-a) \subset (M(-a))_Q$.  Since $M(-a)$
is straight, $M(-a)_Q$ is finitely generated.  Thus $M_\CE$ is finitely
generated, so $M$ has a finitely generated essential submodule by
Proposition~\ref{prop:essential}.  In particular, its Bass numbers in
cohomological degree zero are finite.  Moreover, if $J^\spot$ is a
minimal injective resolution of $M$, then $J^0(-a)$ is straight because
it has a $Q$-graded essential submodule $M(-a)_Q$.  Therefore ${\rm
coker}(M(-a) \to J^0(-a))$ is straight, whence the result follows by
induction on the cohomological degree.  \endproof

\medskip
\noindent
The Bass numbers of such an $M$ at ungraded primes are also finite, by
results of \cite{GWii}.

\begin{thm}\label{thm:finite}
Let $Q$ be an affine semigroup whose saturation is simplicial modulo
units.  If $M$ is a finitely generated $\qgp$-graded $k[Q]$-module, then
the Bass numbers of $H^i_I(M)$ are finite, for any $Q$-graded ideal $I$.
\end{thm}
\proof
This is immediate from Proposition~\ref{prop:strfinite} and Corollary~
\ref{cor:locstraight}.
\endproof

\begin{cor}\label{cor:cechsimplicial}
Suppose $Q$ is affine, with $Q^\sat$ simplicial modulo units, and let $M
\in \ModR$ be finitely generated over $R = k[Q]$.  Then there exist
$\beta \in \qgp$ and $n \in \NN$ such that
$$
  H^i_I(M)(-\beta) \cong \cech \eext^i_{R}(R/I^n, M(-\beta)).
$$
\end{cor}
\proof Since the Bass numbers of $H^i_I(M)$ are finite, by
Remark~\ref{remark:shift} there exists $\beta$ such that
$H^i_I(M)(-\beta)$ is fixed by the \v Cech hull.  Then $H^i_I(M)(-\beta)$
is the Cech hull of its $Q$-graded part and the result follows from
Proposition~\ref{prop:ext}.  \endproof

%\end{section}{Finiteness for simplicial semigroups}%%%%%%%%%%%%%%%%%%%%%
%%%%%%%%%%%%%%%%%%%%%%%%%%%%%%%%%%%%%%%%%%%%%%%%%%%%%%%%%%%%%%%%%%%%%%%%%
%%%%%%%%%%%%%%%%%%%%%%%%%%%%%%%%%%%%%%%%%%%%%%%%%%%%%%%%%%%%%%%%%%%%%%%%%
\section{The essential set}\label{sec:essential}%%%%%%%%%%%%%%%%%%%%%%%%%

The essential set, introduced in Definition~\ref{defn:essential}, fleshes
out in some detail the combinatorics hidden in an affine semigroup.
Since we believe this combinatorics is of independent interest, we
determine in Theorem~\ref{thm:normal}%
, Proposition~\ref{prop:irrelevant}, and the comments in between,
the structure of the essential set in the saturated case, along with its
relation to Hilbert bases%.
, monomial modules, irrelevant ideals of toric varieties, and Alexander
duality.  \comment{; references are needed for each of these}%
The rest of the section we devote to providing the necessary relations
between essential sets of unsaturated semigroups and those of their
saturations, including the results already applied in
Section~\ref{sec:finiteness}.

We need a bit of notation.  Associated to an affine semigroup $Q$ are its
{\em facets} $F_1, \dots, F_r$; these are the degrees of homogeneous
elements outside of the $r$ codimension-one $Q$-graded primes of $k[Q]$.
There are unique primitive integer-valued linear functionals $\{\tau_1,
\dots, \tau_r\}$ on $\qgp$, nonnegative on $Q$, such that $F_i = \{b \in
Q \mid \tau_i(b) = 0\}$.  Given $\alpha \in \qgp$, define $\tau(\alpha)
\in \ZZ^r$ to be the vector $(\tau_1(\alpha), \ldots, \tau_r(\alpha))$,
and let $\tau(\alpha)^+$ be the vector obtained from $\tau(\alpha)$ by
replacing its negative entries with zeros.

\begin{lemma} \label{lemma:tauplus}
Suppose $Q$ is saturated.  Then $(\alpha + Q) \cap Q = (\beta + Q) \cap
Q$ if and only if $\tau(\alpha)^+ = \tau(\beta)^+$.  In particular,
$\varepsilon$ is an essential point if and only if $\tau(\varepsilon)^+
\neq \tau(a + \varepsilon)^+$ for all nonunits $a$ in some generating set
for $Q$.
\end{lemma}
\proof Since $Q$ is saturated, $(\alpha + Q) \cap Q$ is the set of
lattice points $\gamma \in \qgp$ inside the polyhderon defined by
$\{\tau_i(\gamma) \geq 0\ {\rm and}\ \tau_i(\gamma) \geq \tau_i(\alpha)
\mid i = 1,\ldots,r\}$.  This is the polyhedron defined by the
inequalities $\{\tau_i(\gamma) \geq \tau(\alpha)^+_i \mid i =
1,\ldots,r\}$, and the first claim follows easily.

The map $\tau: \qgp \rightarrow \ZZ^r$ takes $Q$ to the semigroup
$\tau(Q)$ isomorphic to the quotient of $Q$ by its group of units.  As a
consequence, the definition of essential point translates to:
$\varepsilon$ is an essential point if and only if $\tau(\varepsilon)^+
\neq \tau(a + \varepsilon)^+$ for all nonunits $a \in Q$.  But since
$\tau_i$ is nonnegative on $Q$ for all $i$, the second statement
follows.~\endproof

\medskip
The lack of nontrivial units in $\tau(Q)$ endows it with a unique minimal
set $\CH$ of semigroup generators, called the {\em Hilbert basis of
$\tau(Q)$}.  Each element of $\CH$ imposes a condition that $\alpha$ must
satisfy to be essential.  To express this condition, define, for $h \in
\NN^r$, the set $\<h\> = \{\zeta \in \ZZ^r \mid \zeta_i > -h_i$ for some
$i$ such that $h_i > 0\}$.  Observe that $\<h\>$ is a union of
half-spaces, and is defined in such a way that $\tau^{-1}(\<\tau(a)\>) =
\{\varepsilon \in \qgp \mid \tau(\varepsilon)^+ \neq
\tau(\varepsilon+a)^+\}$.

\begin{thm} \label{thm:normal}
If $Q$ is saturated then the essential set $\CE$ consists entirely of
essential points.  Furthermore, $\CE = \tau^{-1}(\bigcap_{h \in \CH}
\<h\>) = \bigcap_{h \in \CH} \tau^{-1}(\<h\>)$.
\end{thm}
\proof The second sentence follows from Lemma~\ref{lemma:tauplus} and the
remarks following it because $\tau^{-1}(\CH)$ generates $Q$.  Since
$\<h\>$ is stable under the action of $Q$, $\tau^{-1}(\<h\>)$ is
$Q$-stable, too.  Thus the essential points already form a $Q$-set, which
therefore equals $\CE$.~\endproof

\begin{example} \rm
We consider once again the semigroup $Q$ generated by three elements
$x,y,z$ such that $x + y = 2z$.  $\tau$ embeds $Q$ in $\ZZ^2$ by
sending $x$ to $(2,0)$, $y$ to $(0,2)$, and $z$ to $(1,1)$; the image
of $\qgp$ in $\ZZ^2$ is the sublattice of index $2$ consisting of
pairs $(u,v)$ such that $u$ and $v$ have the same parity.

The points $(0,2)$, $(1,1)$ and $(2,0)$ form a Hilbert basis for
$\tau(Q)$. $\<(0,2)\>$ consists of those points $(u,v)$ in $\ZZ^2$ 
with $v \geq -1$; similarly $\<(2,0)\>$ consists of those points
with $u \geq -1$.  Finally $\<(1,1)\>$ consists of points with
either $u$ or $v$ nonnegative.  Thus the essential points are those
points $\alpha \in \qgp$ such that $\tau(\alpha)$ has one of the following
forms:
\begin{enumerate}
\item $(-1,v)$ for some odd positive $v$
\item $(u,-1)$ for some odd positive $u$
\item $(u,v)$ for $u,v$ nonnegative and with the same parity
\end{enumerate} 

Figure 1 shows the essential set embedded in $\ZZ^2$ via
$\tau$.  The spots (both hollow and solid) represent elements 
of $\qgp$; solid spots are
essential points.  The regions in which $\alpha + Q \cap Q$
remain constant are enclosed by dotted lines.  Note in particular that
there is an essential point in every region, and that the essential
points form a $Q$-set, as predicted by the theorem.

If we refer back to Example~\ref{ex:cech}, we see that the socles
of $R^p\cech$ computed there lie within the essential set, as
Proposition~\ref{prop:essential} and Proposition~\ref{prop:cechstraight}
predict.  Also note that $z + \CE \subset Q$.  One will be able to
translate $\CE$ so that it lies in $Q$ precisely when $Q$ is simplicial;
this is the content of Proposition~\ref{prop:simplicial}, which is the
central goal of this section.  

\begin{figure} 
\includegraphics{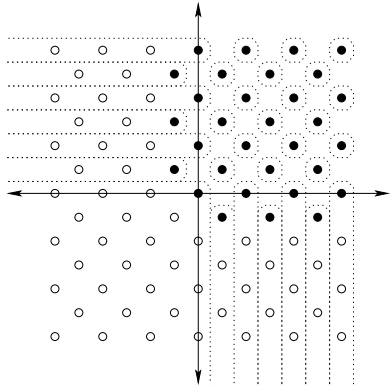}
\caption{The essential set of the semigroup generated by \((2,0)\),
\((1,1)\), and \((0,2)\)}
\end{figure}
\end{example}

%\comment{
\medskip
The essential set $\CE$ is related to a number of other notions already
playing roles in the study of semigroup algebras and toric varieties.
For instance, the subset $\<\CH\> := \bigcap_{h \in \CH} \<h\>$ of
$\ZZ^r$ is a {\em monomial module}%
~\cite{BS},
so $\CE$ might be called a {\em skew monomial module} inside the lattice
$\qgp$.  To get a better picture, $\<\CH\> \subset \ZZ^r$ is a ``fuzzy
neighborhood'' of a certain union $\UU$ of orthants, in the sense that
there is a vector $z \in \NN^r$ such that $\UU \subseteq \<\CH\>
\subseteq \UU - z$.  In fact, each set $\<h\>$ contains and approximates
the union $U_h = \bigcup_{h_i > 0}\{\zeta \in \ZZ^r \mid \zeta_i \geq
0\}$ of half-spaces, and $\UU = \bigcap_{h \in \CH} U_h$; our $z$ can be
any vector with $z_i > h_i$ for all $h\in\CH$ and all $i$.

We can get an even better handle on $\UU$ in the case where $Q$ is the
cone over $\oQ$, an integral polytope in $\ZZ^{d-1} \times \{1\} \subset
\ZZ^{d-1} \times \ZZ = \ZZ^d$.  For each face $F$ of $Q$, let $\van(F)
\subseteq \{1, \ldots, r\}$ be the indices of functionals vanishing on
$F$; similarly, for $z \in \NN^r$, let $\van(z) = \{i \in \{1,\ldots,r\}
\mid z_i = 0\}$.  For instance, if $F = F_i$ is a facet then $\van(F_i) =
\{i\}$, and $\van(z) = \van(F)$ if and only if $z_i = 0$ and $z_j > 0$
for all $j \neq i$.  The polynomial ring $k[\NN^r]$ is the {\em Cox
homogeneous coordinate ring}%
~\cite{Cox}
of the projective toric variety $X$ whose isomorphism class and embedding
in projective space are determined by $\oQ$.  The Cox ring comes equipped
with the {\em irrelevant ideal} $B = \<\xx^z \mid z \in \NN^r\ {\rm and}\
\van(z) \subseteq \van(F) \hbox{ for some face } F\ {\rm of}\ Q\>$.

\begin{prop} \label{prop:irrelevant}
If $\cech$ is the \v Cech hull over $\NN^r$, then $-\zeta \in \UU$ if and
only if $\xx^\zeta \not\in \cech(B)$.  Equivalently, the
$k[\NN^r]$-submodule $\<\xx^\zeta \mid \zeta \in \UU\> \subset k[\ZZ^r]$
is the shift by $(1, \ldots, 1) \in \ZZ^r$ of the \v Cech hull
$\cech(B^\star)$ of the ideal $B^\star$ {\em Alexander dual} to $B$%
~\cite{ER},~\cite[Lecture~VI]{MP}.
\end{prop}
\noindent
In the case where $\ol Q$ is a simple polytope, so the corresponding
projective toric variety is simplicial, $B^\star$ is the Stanley-Reisner
ideal for the simplicial polytope polar to $\ol Q$.

\medskip
\noindent
\proof The equivalence of the two statements is
~\cite[Lemma~2.11]{Mil1}.
Note that $\cech(B)$ is, a priori, a submodule of $k[\ZZ^r]$ since the
latter is the injective hull of $B$ and the former is an essential
extension.  Now $\xx^\zeta \in \cech(B)$ if and only if $\van(\zeta^+)
\subseteq \van(F)$ for some face~$F$ by Remark~\ref{rk:poly}.  On the
other hand, $-\zeta \in U_h$ precisely when $\zeta_i \leq 0$ for some $i$
with $h_i > 0$; that is, when $\van(\zeta^+) \not\subseteq \van(h)$.
Therefore, $-\zeta \not\in \UU$ if and only if $\van(\zeta^+) \subseteq
\van(h)$ for some $h \in \CH$.  This occurs if and only if $\van(\zeta^+)
\subseteq \van(F)$ for some face $F$, because: $\van(F) \subseteq
\van(h)$ for all $h \in F$; and each $h \in \CH$ lies in the relative
interior of some $F$, so $\van(h) = \van(F)$ for this $F$.  We conclude
that $\xx^\zeta \in \cech(B)$ if and only if $-\zeta \not\in
\UU$.~\endproof
%}

\begin{example} \rm
We illustrate this for the semigroup $Q$ of Example~\ref{ex:har};
that is, the semigroup on four generators
$\{x,y,u,v\}$ with the relation $x + u = y + v$.
If we order the four facets appropriately, $\tau$ embeds
$\qgp$ in $\ZZ^4$ by sending $x$ to $(1,0,0,1)$, $y$ to $(0,1,0,1)$,
$u$ to $(0,1,1,0)$ and $v$ to $(1,0,1,0)$.  The image $\tau(\qgp)$ is 
the lattice consisting of points $(a,b,c,d)$ with $a + c = b + d$. 

Now $\<(1,0,0,1)\>$ is the set of $(a,b,c,d)$ in $\ZZ^4$ such
that either $a$ or $d$ is nonnegative.  Thus, if $\alpha \in \qgp$
is an essential point, with $\tau(\alpha) = (a,b,c,d)$, then either
$a$ or $d$ is nonnegative.  Similarly, using the other elements of
the Hilbert basis, we find that:
\begin{itemize}
\item Either $b$ or $d$ is nonnegative.
\item Either $b$ or $c$ is nonnegative.
\item Either $a$ or $c$ is nonnegative.
\end{itemize}

Note that in this example, the irrelevant ideal $B$ of 
Proposition~\ref{prop:irrelevant} is generated
by the elements $\xx^{(1,0,0,1)}, \xx^{(0,1,0,1)}, \xx^{(0,1,1,0)},$ and
$\xx^{(1,0,1,0)}$ in the polynomial ring $k[\NN^4]$.  The {\v C}ech
hull of this ideal is thus supported precisely on those $(a,b,c,d) \in \ZZ^4$
such that at least one of the pairs $\{a,d\}, \{b,d\}, \{b,c\}, \{a,c\}$
consists of strictly positive integers.  Therefore, an element of $\ZZ^4$
fails to be in this support if and only if its negative satisfies the
above four conditions.  To summarize, those elements $\alpha$ of $\qgp$
such that $-\tau(\alpha)$ is not in the support of $\cech B$ are essential
points, as Proposition~\ref{prop:irrelevant} predicts. 

The conditions on $\tau(\alpha)$ given above,
together with the fact that $(a,b,c,d) = \tau(\alpha)$ 
(and therefore $a+b=c+d$), imply that $\alpha$ is an essential 
point if (and only if)
at most one of $\{a,b,c,d\}$ is negative.  Note that no finite shift
will take all of the essential points inside of $Q$, since 
one has essential points $\alpha$ whose negative coordinate is 
$-n$ for any natural number $n$.  In particular, the degrees of the
socle elements of $H^2_{(x,u)}(\omega_R)$ produced in
Example~\ref{ex:har}
are a set of essential points whose negative coordinates are unbounded
below.  

More generally, the fact that we cannot shift the
essential set into $Q$ means that we cannot rule out the
possibility of local cohomology having infinite Bass
numbers.  In fact, we will construct a local cohomology
module with infinite Bass numbers whenever the essential
set cannot be shifted into $Q$
(Corollary~\ref{cor:socle}).
\end{example}

%\begin{problem}
%Find the minimal subsets of $\CH$ yielding $\bigcap_{h \in \CH} \<h\>$
%and $\CE$.
%\end{problem}

\medskip
During our proof of Theorem~\ref{thm:finite}, we needed certain results
about the structure of $\CE$.  In order to obtain them in the generality
we used in Section~\ref{sec:finiteness}, we no longer assume that $Q$ is
saturated.  We begin with a lemma used in the proof of
Proposition~\ref{prop:essential}.
\begin{lemma} \label{lemma:exists}
Given $\alpha \in \qgp$, there exists some essential point $\varepsilon$
with $(\varepsilon+Q) \cap Q = (\alpha+Q) \cap Q$.
\end{lemma}
\proof Any element $\beta \in \qgp$ satisfying $(\beta+Q) \cap Q =
(\alpha+Q) \cap Q$ must also satisfy $\tau(\beta) \preceq \tau(\gamma)$
for all $\gamma \in (\alpha+Q) \cap Q$, where $\preceq$ is the partial
order by componentwise comparison.  The set of possibilities for
$\tau(\beta) \in \ZZ^r$ satisfying this condition is bounded above, and
thus has a maximal element $\tau(\varepsilon)$.  Moreover, $\tau(\beta) =
\tau(\varepsilon) \hbox{ for some } \beta \in \qgp$ if and only if $\beta
- \varepsilon$ is a unit of $Q$.  This proves that $\varepsilon$ is an
essential point.  \endproof

\medskip
The major combinatorial result used in the previous section is the fact
that if $Q^\sat$ is (modulo its units) simplicial, then $\CE$ can be
shifted inside of $Q$.  Therefore, we want an analog to
Lemma~\ref{lemma:tauplus} which holds even for semigroups which are not
saturated.  The key tool relating the combinatorics of a semigroup to the
combinatorics of its normalization is provided by the next lemma.  While
a purely combinatorial proof of this lemma is available, the commutative
algebra below is shorter and more conceptual.  Recall that a {\em face}
of $Q$ is the set of degrees of elements outside a prime ideal of $k[Q]$.
\begin{lemma} \label{lemma:2}
Let $F$ be a face of $Q$.  There exists $a_F \in F$ such that $a_F + Q^F
\subset Q$, where $Q^F := (Q + F^\gp) \cap Q^\sat$ is the {\em partial
saturation of $Q$ at $F$}.
% If $F=Q$ this is just $Q^\sat$; if $F =$ empty face, this is just $Q$.
\end{lemma}
\proof Let $R' = k[Q^F]$ and $\wt R = k[Q^\sat]$.  Then, letting $\pp
\subset R = k[Q]$ be the prime ideal such that $R/\pp = k[F]$, the
$R$-algebra $R'$ is the intersection $R_{(\pp)} \cap \wt R$ of the
homogeneous localization at $\pp$ with the normalization.  The Lemma
calls for a homogeneous element outside of $\pp$ to be in the {\em
conductor ideal}
$$
  {\rm ann}_R(R'/R) = \{x \in R \mid xR' \subset R\},
$$
Such an element exists precisely when ${\rm ann}_R(R'/R)_{(\pp)} = R$;
i.e.\ when the localizations $R_{(\pp)}$ and $R' \otimes_R R_{(\pp)}$ are
equal.  But
$$
\begin{array}{r@{\ =\ }l}
R' \otimes_R R_{(\pp)} 
&	R_{(\pp)} \cap (\wt R \otimes_R R_{(\pp)})\\
&	R_{(\pp)} \cap \wt R_{(\pp \wt R)}\\
&	R_{(\pp)}
\end{array}
$$
because $R_{(\pp)} \subseteq \wt R_{(\pp \wt R)}$.  \endproof

\begin{remark} \rm When $F=Q$, then $Q^F = Q^\sat$, and this is the well-known
fact that every semigroup $Q$ contains an element $a$ with 
$a + Q^\sat \subset Q$.
\end{remark}

\medskip
With Lemma~\ref{lemma:2} in hand, we can now find a sufficient condition
under which $(\beta+Q) \cap Q = (a+\beta+Q) \cap Q$.  Let $\van_\tau(F)
\subseteq \{\tau_1, \ldots, \tau_r\}$ be the subset consisting of
functionals vanishing on the face $F$ of $Q$.

\begin{lemma} \label{lemma:3}
Let $a \in F$, and suppose we have $\beta \in \qgp$ such that
$\tau_i(a+a_F+\beta) \leq 0$ for all $\tau_i \not\in \van_\tau(F)$.  Then
$(\beta+Q) \cap Q = (a+\beta+Q) \cap Q$.
\end{lemma}
\proof Noting that $\tau_j(a) = \tau_j(a_F) = 0$ for $\tau_j \in
\van_\tau(F)$, the hypothesis on $a+a_F+\beta$ implies that the
intersections with $Q$ are contained in $(a+a_F+\beta + Q^\sat)$.
Therefore, it is enough to show that
$$
  (\beta+Q) \cap (a+a_F+\beta\,+\,Q^\sat) = (a+\beta+Q) \cap
  (a+a_F+\beta\,+\, Q^\sat).
$$
This follows by adding $\beta$ or $a+\beta$ to both sides of the equality
in Lemma~\ref{lemma:fact3}, below, and setting respectively $b = a+a_F$
or $b = a_F$.  \endproof

\begin{lemma} \label{lemma:fact3}
If $b \in F$ and $b + Q^F \subseteq Q$, then $Q \cap (b+Q^\sat) = b +
Q^F$.
\end{lemma}
\proof We show $Q \cap (b+Q^\sat) = (b+Q^F) \cap (b+Q^\sat)$, which
obviously equals $b+Q^F$.  Now $Q \cap (b+Q^\sat) \supseteq (b+Q^F) \cap
(b+Q^\sat)$, because $Q$ contains $b+Q^F$; and $Q \cap (b+Q^\sat)
\subseteq (b+Q^F) \cap (b+Q^\sat)$, because $a-b \in Q^F$ when $a \in
\,Q \cap (b+Q^\sat)$, by definition of $Q^F$.~\endproof

\medskip
We are now in a position to state and prove the unsaturated analog of 
Lemma~\ref{lemma:tauplus}.
\begin{prop} \label{prop:tauplusb}
Choose $a_Q$ so that $\tau_i(a_Q) \geq \tau_i(a_F)$ for all $i$ and $F$.
Suppose $a \in Q$, $\beta \in \qgp$, and \hbox{$\tau(a_Q+\beta)^+ =
\tau(a+a_Q+\beta)^+$}.  Then $(\beta+Q) \cap Q = (a+\beta+Q)\cap Q$.
\end{prop}
\proof Let $F$ be the smallest face of $Q$ containing $a$.  Then for all
$\tau_i$ not vanishing on $F$, we have $\tau_i(a) > 0$, so
$\tau_i(a+a_Q+\beta) \leq 0$ (as otherwise the $i^{\rm th}$ coordinates
of $\tau(a_Q+\beta)^+$ and $\tau(a+a_Q+\beta)^+$ are unequal).  Thus
$\tau_i(a+\beta) \leq \tau_i(-a_Q) \leq -\tau_i(a_F)$, and by
Lemma~\ref{lemma:3} we have $(\beta+Q) \cap Q = (a+\beta+Q)\cap Q$, as
required.  \endproof

\medskip
The approximation to Theorem~\ref{thm:normal} in the unsaturated case is
as follows.
\begin{cor}\label{cor:unsaturated}
Let $\CH$ be the Hilbert basis for $\tau(Q)$, and $a_Q$ be as in
Proposition~\ref{prop:tauplusb}.  Then $\CE + a_Q \subseteq
\bigcap_{h \in \CH} \tau^{-1}(\<h\>)$.
\end{cor}
\proof Pick, for each $h \in \CH$, an element $q_h \in Q$ with $\tau(q_h)
= h$.  Suppose $\varepsilon$ is an essential point.  Setting $\beta =
\varepsilon$ and $a = q_h$ in Proposition~\ref{prop:tauplusb}, we have
\hbox{$\tau(\varepsilon+a_Q)^+ \neq \tau(\varepsilon+a_Q+q_h)^+$}.  Just
as before Theorem~\ref{thm:normal}, we have $\tau(\varepsilon+a_Q) \in
\<h\>$, and this holds for all $h \in \CH$.~\endproof

\medskip
\comment{
The following result, which was pivotal in Section~\ref{sec:finiteness},
is a reflection of the fact that the irrelevant ideal in the Cox ring of
a weighted projective space is the maximal ideal.}
\begin{prop}\label{prop:simplicial}
Suppose $Q^\sat$ is simplicial (modulo units).  Then there exists $a
\in Q$ such that $a + \CE \subset Q$.
\end{prop}
\proof The hypothesis on $Q$ means precisely that for each $i =
1,\ldots,r$, the image $\tau(Q)$ contains an element $h^i \in \CH$ in its
Hilbert basis whose unique nonzero coordinate is $h^i_i > 0$.  Observe
that $\<h^i\> = \{\zeta \in \ZZ^r \mid \zeta_i > -h^i_i\}$ is a
half-space by definition.  Setting ${\mathbf h} = (h^1_1, \ldots, h^r_r)
\in \ZZ^r$, we find that ${\mathbf h} + \bigcap_{h \in \CH} \<h\>
\subseteq {\mathbf h} + \bigcap_{i=1}^r \<h^i\> \subseteq \NN^r$.  By
Corollary~\ref{cor:unsaturated} we may take $a = a_Q + \tilde h$ for any
$\tilde h \in Q$ satisfying $\tau(\tilde h) = {\mathbf h}$.
%, where $a_Q$ is as in Proposition~\ref{prop:tauplusb}.
\endproof

%\end{section}{The essential set}%%%%%%%%%%%%%%%%%%%%%%%%%%%%%%%%%%%%%%%%
%%%%%%%%%%%%%%%%%%%%%%%%%%%%%%%%%%%%%%%%%%%%%%%%%%%%%%%%%%%%%%%%%%%%%%%%%
%%%%%%%%%%%%%%%%%%%%%%%%%%%%%%%%%%%%%%%%%%%%%%%%%%%%%%%%%%%%%%%%%%%%%%%%%
\section{Infinite-dimensional socles}\label{sec:infinite}%%%%%%%%%%%%%%%%

In this section we prove our principal result concerning semigroup rings,
Theorem~\ref{thm:converse}, by combining Theorem~\ref{thm:finite} with its
converse, namely that if $Q^\sat$ is not simplicial then one can always
find an ideal $I$ for which local cohomology is not well-behaved.  We
avoid dealing with nontrivial units here, since they add nothing to the
content, but obscure the statement.

\begin{thm} \label{thm:converse}
Let $Q$ be an affine semigroup of dimension $d$ with trivial unit group
(but not necessarily saturated).  The following are equivalent:
\begin{enumerate}
\item
$Q^\sat$ is simplicial.

\item
For every $Q$-graded ideal $I$ and every finitely generated $Q$-graded
$k[Q]$-module $M$, the Bass numbers of $H^i_I(M)$ are finite.

\item
For every $Q$-graded prime $\pp$ of dimension~2,
$H^{d-1}_\pp(\omega_{k[Q^\sat]})$ has finitely generated socle.
\end{enumerate}
\end{thm}

This theorem provides a proof and generalization of Example~\ref{ex:har}.
The key to our argument is Yanagawa's computation of the local cohomology
of the canonical module $\omega_{k[Q]}$ over a normal semigroup ring
$k[Q]$ \cite{YanPoset}.  To state it, let $\tau_1, \dots, \tau_r$ be
linear functionals which vanish on the facets $F_1, \dots, F_r$ of $Q$
and take nonnegative integer values on $Q$, as in the previous sections.
For the sake of simplicity we assume that $Q$ has no nonzero units.
Choose a hyperplane $H$ transverse to the real cone $\RR_+ Q$ generated
by $Q$, so that $\oQ = (\RR_+ Q) \cap H$ is a polytope of dimension $d-1
= \dim(k[Q]) - 1$ whose faces (including the empty face $\nothing$)
correspond to the primes of $k[Q]$.

\begin{defn} \sl
Let $F \in Q$ correspond to $\oF \in \oQ$ (so $\0 \in Q$ corresponds to
$\nothing \in \oQ$, for example).  Define the polyhedral cell subcomplex
%$$
%  \oQ(\alpha) = \{\ol{F'} \in \oQ \mid (\alpha + \RR_+ Q) \cap F' =
%  \nothing\}
%$$
%and
$$
  \oF(\alpha) = \{\ol{F'} \in \oF \mid (\alpha + \RR_+ Q) \cap F' =
  \nothing\}% = \oF \cap \bigcup\{\oF_i \mid \tau_i(\alpha) > 0\}
$$
of $\ol F$ for any face $\ol F \subseteq \ol Q$ and $\alpha \in \qgp$.
\end{defn}

\begin{thm}\label{thm:yancanonical} \cite[$\!$Theorem 6.1]{YanPoset}
Let $Q$ be saturated and $\pp$ be a graded prime of $k[Q]$, corresponding
to a face $\oF$ of $\oQ$.  Then $H^{d-i}_\pp(\omega_{k[Q]})_\alpha \cong
\HH_{i-1}(\oF,\oF(\alpha))$ for all $\alpha \in \qgp$.~\endproof
\end{thm}

We will apply this when $\pp$ has dimension~$2$; that is, when $\oF$ is
an edge of $\oQ$.  
%In this case the complexes in question are easy to
%compute.

\begin{prop}\label{prop:twoface}
If $\pp$ corresponds to the edge $\oF$ and $Q$ is saturated, then:
\begin{enumerate}
\item
$H^{d-1}_\pp(\omega_{k[Q]})_{\alpha} = 0$ if $\tau_i(\alpha) > 0$ for
some $i$ such that $\ol{F}_i \cap \oF \neq \nothing$.

\item
$H^{d-1}_\pp(\omega_{k[Q]})_{\alpha} = 0$ if $\tau_i(\alpha) \leq 0$ for
all $i$ such that $\ol{F}_i \cap \oF = \nothing$.

\item
$H^{d-1}_\pp(\omega_{k[Q]})_{\alpha} = k$ if neither of the above
conditions holds.
\end{enumerate}
\end{prop}

\proof Suppose the first condition holds.  If $\ol{F}_i$ contains $\oF$,
then $\alpha + \RR_+ Q$ misses $F$ entirely, so $\oF(\alpha) = \oF$, and
the zeroth relative homology is zero.  Otherwise, $\ol{F}_i \cap \oF$ is
a vertex of the edge $\oF$, and $\oF(\alpha)$ contains at least that
vertex.  Thus the relative homology is again zero.

If the second condition holds, but the first does not, then
$\tau_i(\alpha) \leq 0$ for all~$i$.  This implies $\oF(\alpha)$ is the
void complex---not even $\nothing \in \oF(\alpha)$, so the zeroth
relative homology is still zero.

In the third case, $\oF(\alpha)$ consists of just the empty face
$\nothing$, and the zeroth relative homology is the number of connected
components of $\oF$.  \endproof

\begin{cor}\label{cor:socle}
Let $\oF$ be an edge of $\oQ$ such that there exists a facet $\ol{F}_j$
of $\oQ$ with $\ol{F}_j \cap \oF = \nothing$.  Let $\pp$ be the prime of
$k[Q]$ corresponding to $\oF$.  Then if~$Q$ is saturated,
$H^{d-1}_\pp(\omega_{k[Q]})$ has an infinite-dimensional socle.
\end{cor}
\proof Every nonzero element $x \in H^{d-1}_\pp(\omega_{k[Q]})$ is
annihilated by some power of the maximal ideal of $k[Q]$.  To see this,
suppose $x$ is homogeneous of degree $\alpha$, and assume that for some
$\beta \in Q$, we had $\xx^{n\beta}x \neq 0$ for all $n \in \NN$.  Then,
by Proposition~\ref{prop:twoface}, $\tau_i(n\beta + \alpha) \leq 0$ for
all~$n$ and all~$i$ such that $\ol{F}_i \cap \oF \neq \nothing$.  Thus
$\tau_i(\beta) = 0$ for all such $i$, so $\beta \in \ol{F}_i$ for all
such $i$.  But the intersection of all such $\ol{F}_i$ is empty, so
$\beta = \0$.

If $\ol{F}_j \cap \oF = \nothing$, then choose $\alpha \in \qgp$ such
that $\tau_j(\alpha) > 0$ and $\tau_i(\alpha) \leq 0$ for $i \neq j$.
Proposition~\ref{prop:twoface} implies that
$H^{d-1}_\pp(\omega_{k[Q]})_{\alpha}$ is nonzero and killed by some power
of the maximal ideal, so $H^{d-1}_\pp(\omega_{k[Q]})$ has nontrivial
socle.  Suppose its socle were finite-dimensional.  Then there would
exist $\beta \in \qgp$ such that $\tau_j(\beta)$ is maximal among the
socle degrees in $\qgp$.  But since $\tau_j(\alpha) > 0$, we have
$\tau_j(\beta) > 0$, so $\tau_j(2\beta) > \tau_j(\beta)$.  Moreover the
local cohomology is nontrivial in degree $2\beta$.  Taking a nonzero
element of $H^{d-1}_\pp(\omega_{k[Q]})_{2\beta}$ and multiplying it by a
sufficiently large power of the maximal ideal then yields a socle element
in a degree $\gamma$ with $\tau_j(\gamma) > \tau_j(\beta)$, which is a
contradiction.  \endproof

\bigskip
\noindent
{\it Proof of Theorem~\ref{thm:converse}}. $1 \implies 2$ is
Theorem~\ref{thm:finite}, and $2 \implies 3$ is because
$\omega_{k[Q^\sat]}$ is finitely generated over $k[Q]$.  For $3 \implies
1$, the unsaturated case follows from the saturated case.  Indeed, any
$Q^\sat$-graded ideal $I \subset k[Q^\sat]$ is generated up to radical by
elements $\yy = (y_1, \ldots, y_s)$ in $k[Q]$ (high powers of any
homogeneous generating set for $I$ will do).  If $M$ is any
$k[Q^\sat]$-module, the cohomology of the \v Cech complex $C^\spot(\yy;
M)$ on these generators is therefore a module over both $k[Q^\sat]$ and
$k[Q]$.  As such, it is simultaneously the local cohomology of $M$ over
$k[Q^\sat]$ with support on $I \subset k[Q^\sat]$ and over $k[Q]$ with
support on $I \cap k[Q]$.  Furthermore, any socle element of a
$k[Q^\sat]$-module is also a socle element over $k[Q]$, since the maximal
ideal of $k[Q]$ is contained in the maximal ideal of $k[Q^\sat]$.

Thus by Corollary~\ref{cor:socle}, it suffices to produce, for any
polytope $\oQ \neq$ simplex, an edge $\oF$ of $\oQ$ that misses some
facet.  Equivalently, it suffices to show that if $\oQ$ is a polytope in
which every edge meets every facet then $\oQ$ is a simplex.  Let $\oF \in
\oQ$ be a facet, and $\tau$ a linear functional supporting $\oF$,
nonnegative on $\oQ$.  Suppose $\tau$ takes a minimal nonzero value at a
vertex $v \not\in \oF$.  If more than one vertex of $\oQ$ lies off of
$\oF$, there is an edge (necesasarily missing $\oF$) connecting $v$ to
some vertex at which $\tau > 0$.  Thus, if every edge meets every facet,
there can be only one vertex of $\oQ$ lying off of each facet, and $\oQ$
must be a simplex.~\endproof

%\end{section}{Infinite-dimensional socles}%%%%%%%%%%%%%%%%%%%%%%%%%%%%%%
%%%%%%%%%%%%%%%%%%%%%%%%%%%%%%%%%%%%%%%%%%%%%%%%%%%%%%%%%%%%%%%%%%%%%%%%%
%%%%%%%%%%%%%%%%%%%%%%%%%%%%%%%%%%%%%%%%%%%%%%%%%%%%%%%%%%%%%%%%%%%%%%%%%
\section{Open problems}\label{sec:problems}%%%%%%%%%%%%%%%%%%%%%%%%%%%%%%

It has been seen above that affine semigroup rings provide a wealth of
examples and counterexamples to general questions about local cohomology
in singular varieties.  In particular, they shed some light on some of the 
general questions posed by Huneke on local cohomology \cite{Hun}:
\begin{enumerate}
\item When is $H^i_I(M)$ zero?
\item When is $H^i_I(M)$ finitely generated?
\item When is $H^i_I(M)$ artinian?
\item When is the number of associated primes of $H^i_I(M)$ finite?
\end{enumerate}
Although the answer to the fourth is trivially ``always'' in the cases
discussed in this paper, the above examples provide clues as to how to
refine the first three, given a grading.

Section~\ref{sec:semigroups} provides a possibility for answering
question~1:\ relate the vanishing of local cohomology in a given
cohomological degree to the vanishing of $\ext$ modules in that
cohomological degree and lower.  Theorem~\ref{thm:extspectral}
establishes this link for graded modules over semigroup rings; we believe
that such a connection exists in significantly more generality, but we
are unaware of how to relate infinitely generated modules to finitely
generated ones without resorting to a grading.  The key concept is that
of a certain kind of ``constancy'', provided here by the Cech hull.  This
type of constancy is reminiscent of the characteristic~0 regular local
case, in which the modules in question are treated as $D$-modules
\cite{Lyu1}.  Perhaps the right generalization of $D$-module to the
singular setting will provide the appropriate notion of constancy to
bridge finitely generated Ext modules and local cohomology.

A partial answer to question~2 in the general (local, ungraded) case
concerns numerical criteria on the heights of primes and cohomological
degress involved \cite{Hun}.  In the semigroup-graded case, finite
generation can be viewed as a convex-geometric problem, dealing with
$Q$-graded degrees in which the summands in a minimal injective
resolution are nonzero.  We expect in the $Q$-graded case for these
considerations to yield geometric and combinatorial criteria in addition
to the general numerical criteria.  For the canonical module of a normal
semigroup ring, for instance, local cohomology at a graded prime ideal
$\pp$ of $R$ is finitely generated if and only if it vanishes, since
Proposition~\ref{prop:canonical} expresses such cohomology in terms of
derived functors of $\cech$, which are never finitely generated if they
are nonzero, or in terms of $\cech \omega_{R/\pp}$, which is also never
finitely generated.

As pointed out by Huneke \cite{Hun}, question 3 has two parts, namely:

3a. When is the maximal ideal the only associated prime of $H^i_I(M)$?

3b. When are the Bass numbers of $H^i_I(M)$ finite?

\noindent
Both 3a and 3b should have concrete combinatorial answers in the
semigroup case, at least when $M$ is a canonical module.  In fact, we
expect the essential set to play a pivotal role in answering these and
the following refinement of 3a: For which cohomological degrees $i$ and
graded ideals $I$ is a given prime of $k[Q]$ associated to
$H^i_I(\omega_{k[Q]})$?

As for question~3b, it seems to be connected with the kinds of
singularities which appear in the normalization of the ring $R$.  Whether
this holds in more generality than simply for semigroup rings is an
interesting question.  For instance, one can try classifying the
singularities of a ring $R$ or ideals $I$ for which the modules
$H^i_I(M)$ can have infinite Bass numbers (or, for that matter, which
primes can appear with infinite Bass number).  Even in the case of a
semigroup ring, we do not have satisfactory answers to these last
questions.

% There is also the question of generalizing our techniques to more
% general semigroup-graded algebras.  The formalism of the \v Cech hull
% works in complete generality in such cases, but it is unclear what of
% the remainder of our theory applies.  What is the proper analogue of
% straight modules in such a setting?  What properties do they have?

Finally, is there a global version of the \v Cech hull that works for
toric varieties, and if so, what is its relation to the \v Cech hull over
the Cox homogeneous coordinate ring \cite{Cox}?  More generally, for
varieties with a torus action, can a global \v Cech hull give information
about cohomology with support on subvarieties fixed pointwise by
subgroups of the acting torus?  In the toric case, properties of a global
\v Cech hull will be governed by the group of Weil divisors modulo
Cartier divisors, introduced by Thompson to control resolutions of
singularities \cite{Tho}.

%\begin{example} \rm
%If $\pp \subseteq k[Q]$ is prime of codimension $c \leq n = \dim(k[Q])$
%then if $i > c$ and some other prime $\qq$ is associated to
%$H^i_\pp(\omega_{k[Q]})$ then the Bass number of this local cohomology
%module at $\pp$ is infinite.  \comment{proof?}
%\end{example}

%\bigskip
%\begin{textlist}
%\item

%\item
%A word on how to globalize for toric varieties, or more general varieties
%with torus actions, if it's just as easy.  Mention smooth toric variety
%case = \v Cech hull on Cox ring (if true).
%\end{textlist}

\bigskip
\noindent
\textsc{Acknowledgements.} The authors wish to thank Bernd Sturmfels and
Kohji Yanagawa for their useful suggestions.  Both authors were partially
supported by the National Science Foundation, and the second was
partially supported by the Alfred P. Sloan Foundation.

%\end{section}{Open problems}%%%%%%%%%%%%%%%%%%%%%%%%%%%%%%%%%%%%%%%%%%%%
%%%%%%%%%%%%%%%%%%%%%%%%%%%%%%%%%%%%%%%%%%%%%%%%%%%%%%%%%%%%%%%%%%%%%%%%%
%%%%%%%%%%%%%%%%%%%%%%%%%%%%%%%%%%%%%%%%%%%%%%%%%%%%%%%%%%%%%%%%%%%%%%%%%
%\bibliographystyle{amsalpha}
%\bibliography{biblio}
%\nocite{YanPoset}
\providecommand{\bysame}{\leavevmode\hbox to3em{\hrulefill}\thinspace}

%%%%%%%%%%%%%%%%%%%%%%%%%%%%%%%%%%%%%%%%%%%%%%%%%%%%%%%%%%%%%%%%%%%%%%%%%
%%%%%%%%%%%%%%%%%%%%%%%%%%%%%%%%%%%%%%%%%%%%%%%%%%%%%%%%%%%%%%%%%%%%%%%%%
\end{document}